S. I. Kryuchkov

# The Four Color Theorem and Trees

Abstract

Connection of the Four Color Theorem (FCT) with some operations on trees is described. L.H.Kauffman's theorem about FCT and vector cross product is discussed. Operation of transplantation on trees linked with the move of brackets according to the associative law is used to formulate a conjecture. When map is represented as a tying of the trees this conjecture proposes the existence of special coloring of this map. This coloring makes possible successive transplantations such that one of these trees is transformed to another, and all the intermediate maps are colored properly. It means (in terms of L.H.Kauffman) that not only for two different positions of brackets in a product of *n* factors the associative law works (on special values of factors) but also there is a way of moving brackets that for all intermediate positioning of brackets the associative law is obeyed.
  Some classes of trees for which the conjecture is proved are presented.

The four color theorem (FCT) has a long history and a lot of equivalent formulations [1-3]. The complexity of its computer proof [4,5] provokes attempts of proving it in a more human way. Now they are mainly new formulations of FCT (see for example [6]). One of such new formulations [7] attracted our attention by its connection with the associative law. The importance of generalized the associative law [8] in various fields is now obvious (see for example [9-11]).

We shall first introduce some facts to establish terminology. It is sufficient to prove the FCT for planar bridgeless maps with strictly 3 countries incident to every vortex [1]. Further on we shall treat only such maps without special mention. For such maps the FCT is equivalent to the problem of coloring edges in 3 colors in such a way that all (3) edges incident to common vortex get different colors (which is called proper coloring [1.3].

Let $E,I,J,K$ be four colors of countries. We shall assign colors to edges according to the rules (edges get the color of the result of multiplication of colors of adjacent countries): $EE = E$, $EI = I$, $EJ = J$, $EK = K$, $II = JJ = KK = E$, $IJ = K$. Those are the multiplication rules of Klein's four-group. It is sufficient to solve the FCT for maps, having the form of two tied trees [3, 7]. See [7] and fig.1 giving clear understanding of what we mean by that.

Every tree is clearly associated with the bracketed product as it is seen in fig1. We use the * operation to get mirror image of bracketed product: $[((xy)z)]* = z(yx)$, and $T$ operation to get tree from such product, and # to mark the tying of two trees. All this is illustrated by fig.1 (see also [7]).

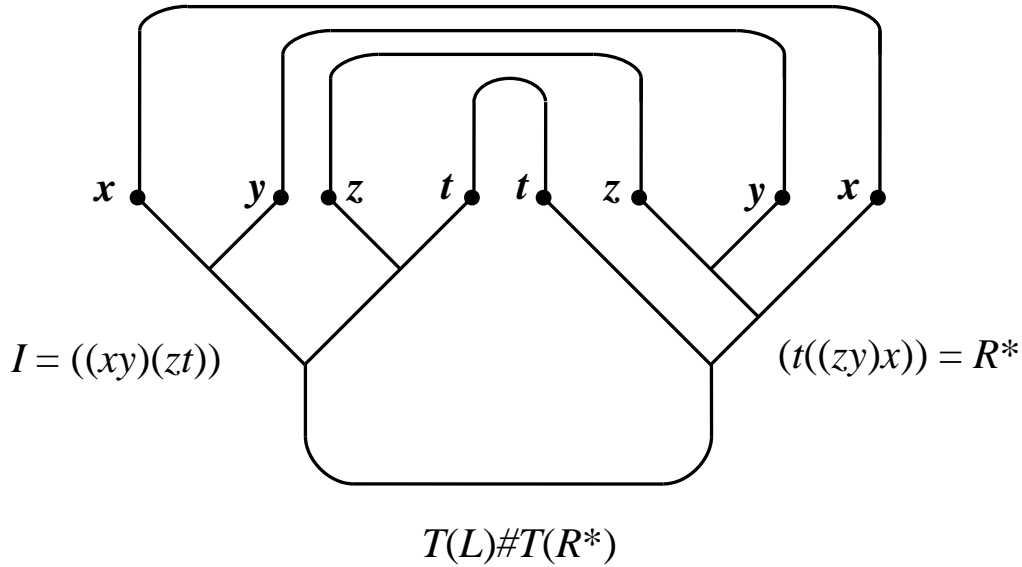

Fig.1. Tying trees.

Let $L$ and $R$ be some bracketed products of variables $x_1, \ldots, x_n$. We say (after L.H.Kauffman [7]) that the solution of $L = R$ is sharp if both sides are nonzero, variables $x_1, \ldots, x_n$ are chosen from $\{i, j, k\}$ – orts in 3-dimentional space, and multiplication is the cross product with the multiplication rules: $ii = jj = kk = 0$, $ij = k, ji = -k$ (and circular permutations of $i, j, k$).

With every equation $L = R$ is associated the $T(L)\#T(R*)$ map. This fact gives rise to the theorem of L.H.Kauffman [7] that FCT is equivalent to the existence of sharp solution of equation $L = R$ for any $n$ and all choice of $L$ and $R$.

Turning to edge coloring problem one can see that when the map is colored properly the edges of map $T(L)\#T(R^*)$ corresponding to $x_1, \ldots, x_n$. in $L$ and $R$ get the colors $I, J, K$. Attributing values $i, j, k$ to $x_i$ when corresponding edges have colors $I, J, K$ we get the sharp solution of $L = +R$ or $L = -R$ (which is clear from comparison of multiplication rules of Klein's four group and those of gross product). In order to prove that only the first case takes place L.H.Kauffman used the notion of formation (see part III of [7]).

This result can be proven in a more short way. It is clear that for sharp solution (and proper coloring) products $ii = jj = kk = 0$ are not encountered anywhere. So only the second group of relations ($ij = k, ji = -k, \ldots$) is used. They are the same as for the Hamiltonian quaternions. So we can treat $L$ and $R$ as quaternion products in this case. Quaternion algebra is associative and $L = R$ is valid on quaternions. Thus the question with the sign in $L = \pm R$ is solved.

The sharp solution of $L = R$ on quaternions can be treated as the solution with only $\{\pm i, \pm j, \pm k\}$ values appearing on every step of calculation according to brackets of $L$ and $R$. Associativity is lost here because the quaternion unit element $e$ is forbidden in use. One is always to be in vector subspace of the 4-dimentional space of quaternion algebra (which is associative). Just the same way we can introduce the projection operator on the set $\{I, J, K\}$. If one associates this operator with the closing bracket then the sharp solution of equation $L = R$ can be introduced for Klein's four group (the solution when $E$ never appears in calculations according to brackets). Klein's four-group is more convenient object because of the commutativity of multiplication. Of course, the equivalence of such sharp solution to FCT is still true.

Thus the FCT gives positive answer to the question, whether two different positioning of brackets in accordance with the associative law may hold the nonzero value of product $(x_1 \ldots x_n)$. (For proper values of $x_i$).

Elementary move of brackets according to the associative law $(A(BC)) <-> ((AB)C)$ corresponds to transformation of the tree. We shall call it transplantation [9]. In [9] transplantation in the tree describing angular momentum addition is associated with $6j$-symbol. Transformation of one tree to another is associated with $3nj$-symbol (which is the sum of products of $6j$-symbols). So this may lead to some new "physical" formulation of FCT which becomes a theorem of quantum theory of angular momentum (or conformal field theory? see [10]). If product contains many $x_i$ then whole sub-tree can be transplanted when two brackets are moved according to the associative law (see fig.2).

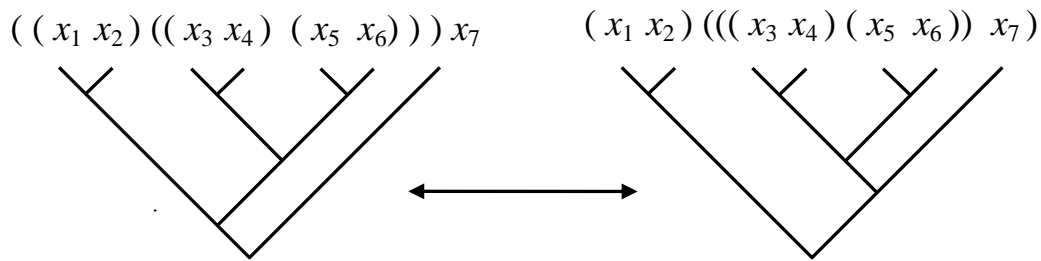

Fig.2. Transplantation of a sub-tree.

When $T(L)$ and $T(R)$ differ only by one transplantation and the chosen values of $x_i$ give sharp solution of $L = R$ we call such transplantation <u>admissible</u>. Sometimes a number of successive transplantations are admissible for chosen values of $x_i$ they form <u>admissible path</u>. In treating the generalized associative law [8] the 5-cucle, or pentagon, diagram is of large importance. In fig.3 adjacent graphs differ only in one transplantation. If such diagram is commutative (for any values of $x_i$) the associative law holds for any $n$ in $(x_1 \ldots x_n)$ [8].

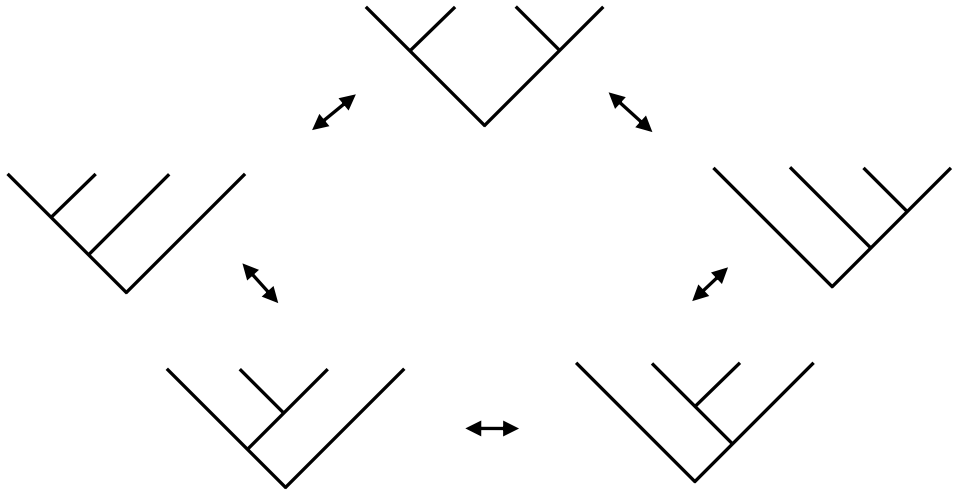

Fig.3. 5-cycle of transplantations.

In our case this is not true. But one can easily see that not only do exist the sharp solutions of $L = R$ for all pairs of trees in fig.3, but also for every such pair exists the succession of admissible transplantations, transforming the tree $T(L)$ into $T(R)$. This admissible path is always the shortest path, and there are no sharp solutions of $L = R$ which make not the shortest path admissible.

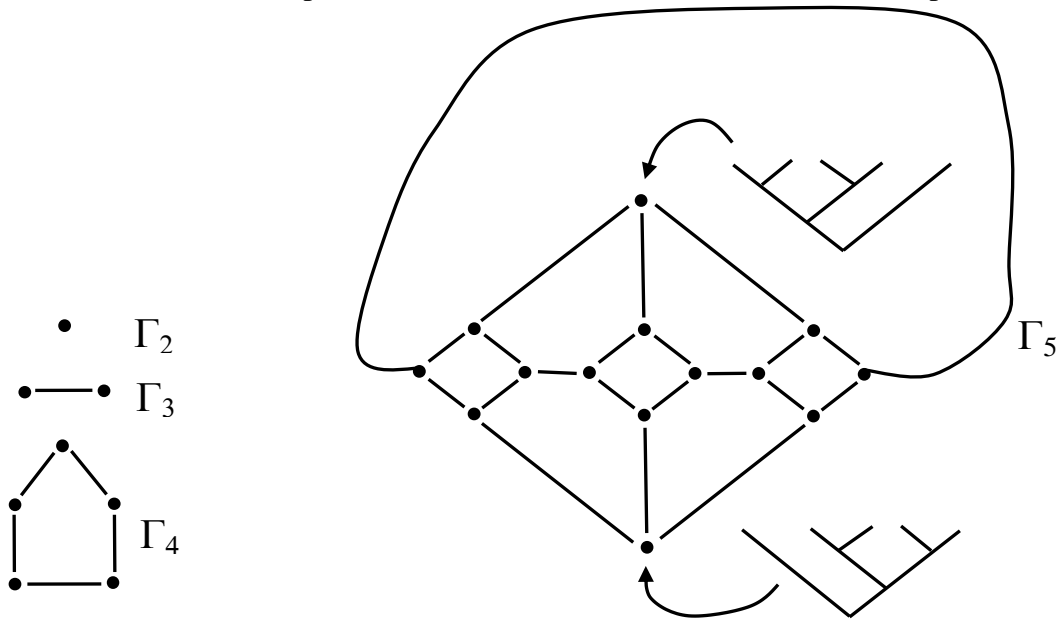

Fig.4. Graphs $\Gamma_n$ for $n \leq 5$.

It is convenient to introduce graph $\Gamma_n$ the vortices of which are all possible trees associated with bracketed products of $n$ factors $x_i$ and edges correspond to transplantations. Figure 4 shows $\Gamma_n$ for $n \leq 5$ and two trees associated with corresponding vortices of $\Gamma_5$. The shortest cycles in $\Gamma_n$ are 5 - and 4 - cycles (see [8, 9, 10]).

Now admissible path becomes <u>admissible pass</u> on graph $\Gamma_n$. It is interesting to know, whether admissible paths exist for any two vortices of $\Gamma_n$. For $n \leq 5$ the shortest paths on $\Gamma_n$ are admissible (of course, for proper choice of $x_i$) which is possible to check directly.

Some steps of path on $\Gamma_n$ commute. This commutativity is associated with quadrangles of $\Gamma_n$ (they correspond to transplantations which do not influence each other). If to exclude the waste steps occurring on some path and corresponding to passing opposite sides of quadrangle, then for $n = 5$ the rule (which is correct for $n = 3, 4$) that only the shortest paths are admissible is true. But for $n > 5$ this rule fails. Though, locally, when transplantations of admissible path are made with no more than 5 sub-trees this shortest-path rule is obeyed.

It is easy to produce admissible paths transforming any graph to the graph of special type as in fig.5. It corresponds to left most position of brackets: $(\ldots((x_1x_2)x_3)\ldots x_n)$. This path consists of transplantations made with the most left sub-trees that can be transplanted to the left. Fig.6 shows an example of such path.

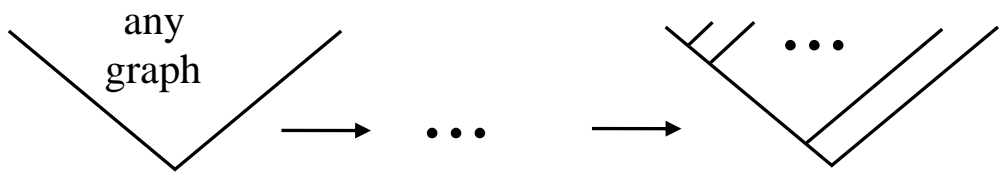

Fig.5. The first example when existence of admissible path is proved.

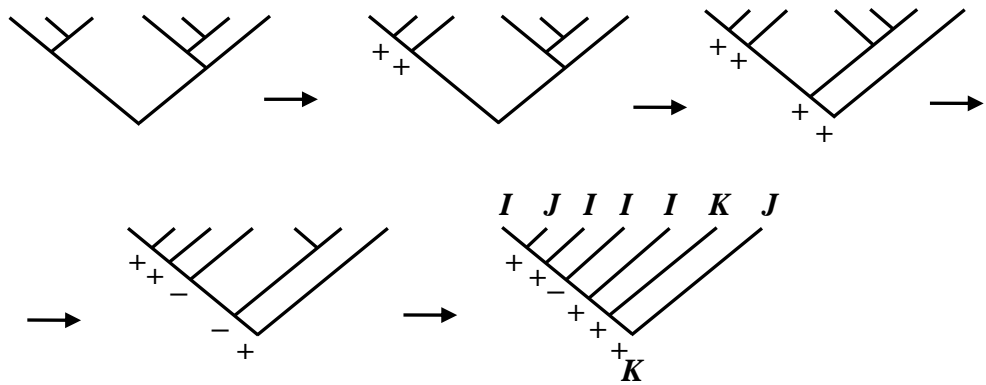

Fig.6. Example of coloring a tree at every step of a path on $\Gamma_n$.

It is convenient to use not the coloring of edges, but colorings of vortices of tree. For sharp solution there are two possibilities for orderings of colors in a vortex: $(IJ) = K$ or $(JI) = K$ (and circular permutations). Correspondingly, we mark vortices by "+" and "−". The rule of treating colors of vortices under transplantations is clear from fig.7. Attributing of signs to vortices is connected with one of the Heawood's formulations of the FCT (see $C_{13}$ in [7]).

admissible transplantations:

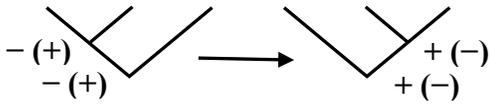

not admissible transplantations:

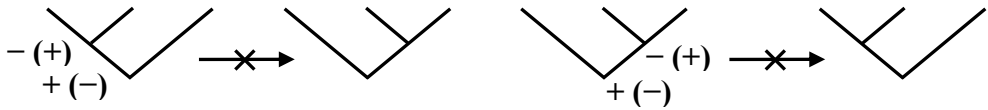

Fig.7. Transplantations and coloring of vortices.

In order to get coloring coherent with transplantations one can assign colors (+ or −) to vortices at every step of the chosen path. For paths of fig.5 type the color assigned at every step can be either arbitrary (if vortices participating in transplantation are not colored) or coordinated with the already colored vortex (as in fig.7). Before every step of fig.5 at least one vortex participating in it is not colored. So the proper coloring is always possible. One of such colorings is shown in fig.6.

Of course, the same procedure is proper for pare of trees, one of which corresponds to rightmost position of all brackets. Also this is true when $L$ or $R$ is the product of rightmost bracketed block on leftmost bracketed block: $(x_1(\ldots(x_{k-1}x_k)\ldots))((\ldots(x_{k+1}x_{k+2})\ldots)x_n)$. This is clear from the correspondence of trees and polygons (see, for example, [12], ch.20), which is shown in fig.8. $x_1$, $x_2$, …,$x_n$ and their product form sides of polygon ($n+1$ - gone). Diagonals are placed in such a way that any two factors and their product form the triangle. All the sides and diagonals of polygon are attributed colors $I$, $J$, $K$ and triangles get signs + and − coherently with the coloring of tree. Transplantation is now the move of one diagonal inside the corresponding quadrangle and it is admissible if the triangles incident to diagonal are of the same sign. The two newly formed triangles get the changed sign.

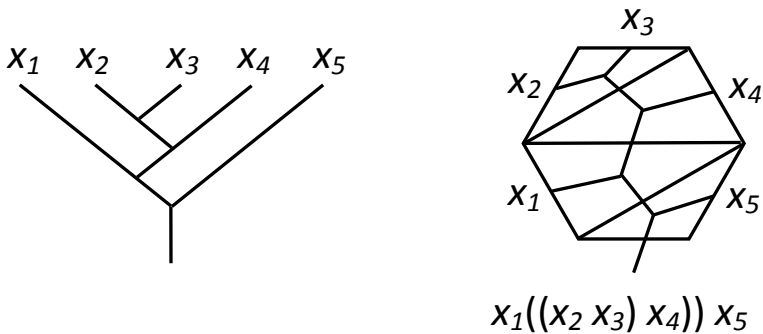

Fig.8. Trees and polygons.

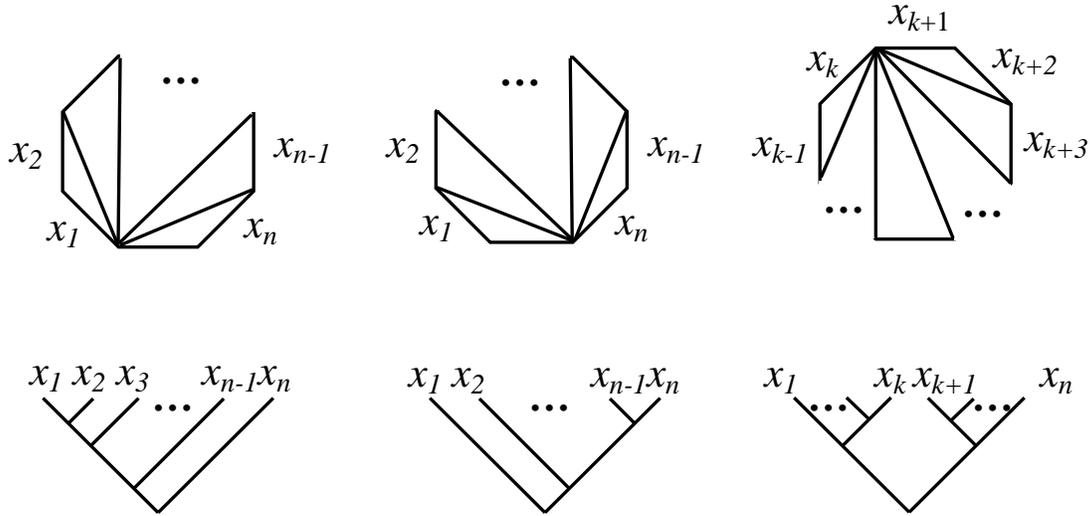

Fig.9. Special class of trees.

The products with the leftmost, rightmost positions of brackets, and their combination, as described above, correspond to trees and polygons of figure 9.

The polygonal formulation is more symmetric one, as far as it is not important which side of polygon is taken as a root of the tree. Polygons a,b,c in fig.9 are congruent and differ only on the position of the vortex incident to all the diagonals. It is clear, that admissible transformation of polygon with arbitrary position of diagonals to any one from fig.9 is carried out by same algorithm. For "a" it was already described. For case "b" it is enough to change "left" for "right" in the algorithm for 'a". It is enough to rotate polygon 'c" (and the one to be transformed) up to coincidence with polygon "a" or 'b", make all transformations according to algorithm and rotate it backwards.

Experimenting with trees and admissible paths on $\Gamma_n$ brought us to a

Conjecture… For any $n \in \mathbb{N}$, $n \geq 3$ equation $L = R$ has such sharp solution that $T(L)$ and $T(R)$ are connected with admissible path on $\Gamma_n$ for this solution.

This conjecture entails FCT but is not equivalent to it. To illustrate this we shall prove the following

Proposition. There exist $L$, $R$ and such sharp solution of $L = R$, that both in $T(L)$ and $T(R)$ not a single transplantation is admissible for this solution.

The proof. First assign + or − signs to every vortex of every tree so that adjacent vortices get opposite signs and the root vortex (corresponding to the last multiplication in $(x_1…x_n)$) is marked (say) by +. Assign to the root edge value (say) $K$. This fixes one of the proper colorings for every tree. In the tree with such coloring no transplantation is admissible.

The number of trees associated with $(x_1…x_n))$ is [12]:

$$g(n) = \frac{(2n-2)!}{(n-1)!n!} \sim \frac{4^{n-1}}{\sqrt{n(n-1)}\, n}$$

Every $x_i$ equals $I$, $J$ or $K$. So there are no more than $3^n$ different attributing of values to all $x_i$. When $g(n) > 3^n$ (such $n$ exists, which is clear from asymptotic of $g(n)$) there are at least two different trees with the same set of $x_i$. Take one of them for $T(L)$, the other for $T(R)$. The set of $x_i$ values gives the sharp solution of $L = R$. Thus the proposition is proved.

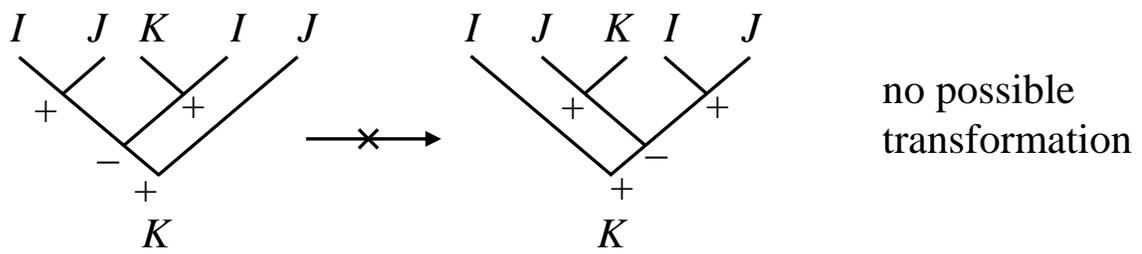

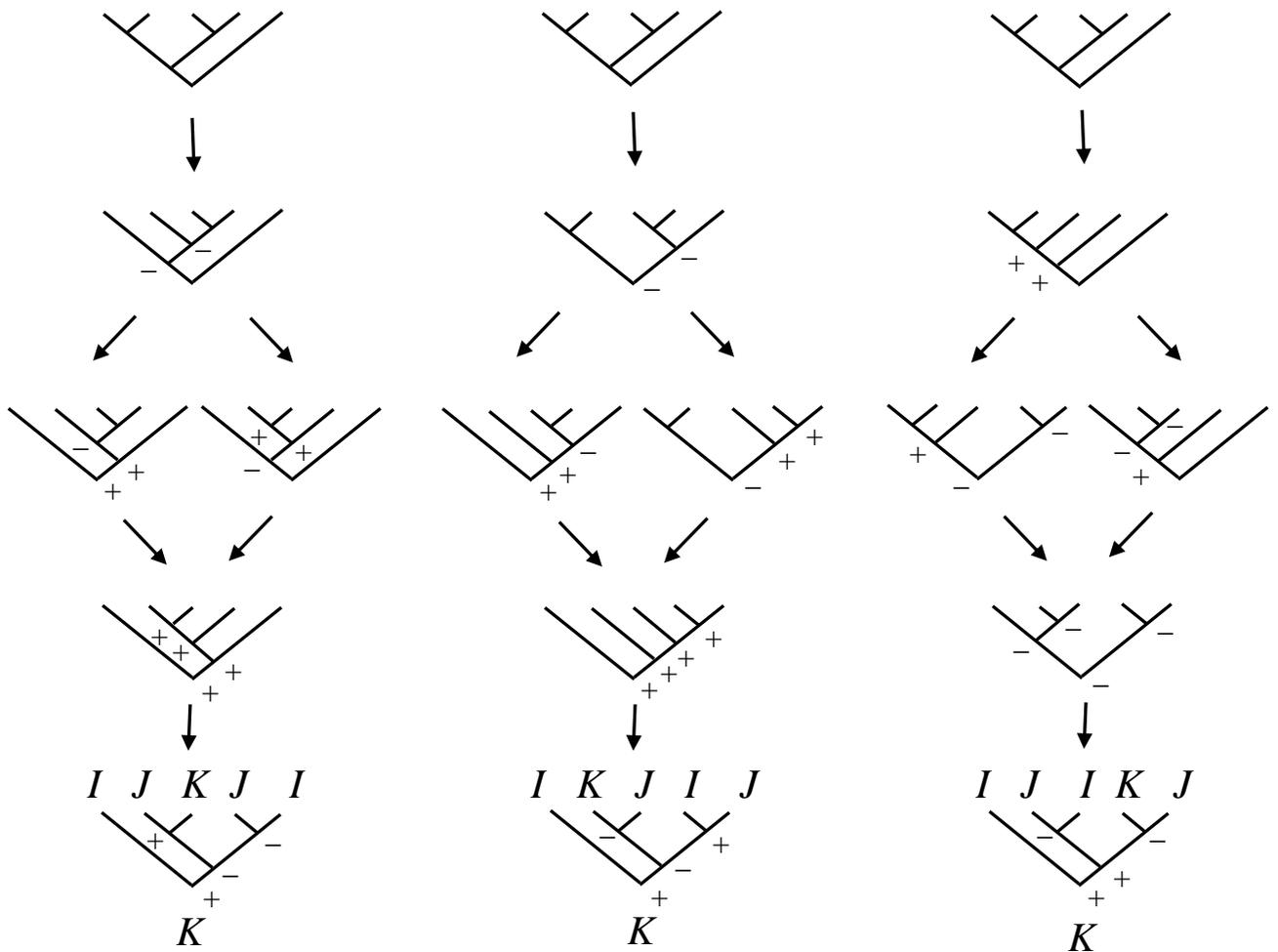

Fig.10. Example for proposition.

One can easily find an example of trees existence of which is stated by proposition. It is clear from the proof that the number of such trees is growing with *n* very strong. The minimal value of n for such tree is overestimated. Figure 10 gives the described example. It also shows that for the same pair of trees more than one path on $\Gamma_n$ (described in conjecture) may exist. The process of coloring vortices is also illustrated in this figure.

When some diagonals of polygons corresponding to *T(L)* and *T(R)* are the same, the problem of finding admissible path factorizes. In this case there are two sub-trees in every tree and transplantations are to be made only inside those sub-trees. So the problem for *n* factors $x_i$ is reduced to problem for trees with the number of $x_i$ less than $n$. This situation is illustrated by fig.11.

We shall present one more class of trees to which any tree can be transformed by admissible path. First we note that the coloring of one (any) vortex of a tree can be taken arbitrary, as well as the coloring of one of the edges of the tree.

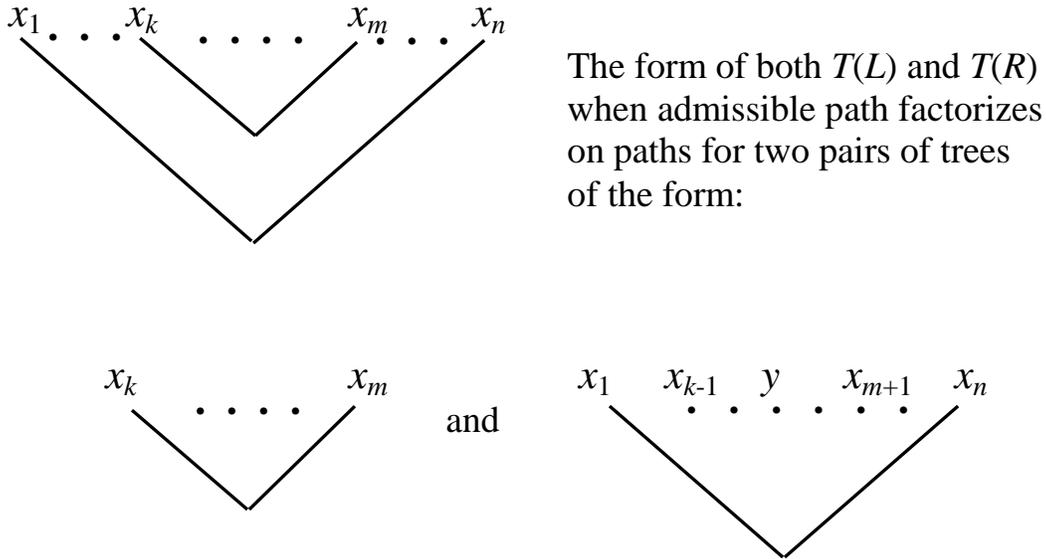

The form of both $T(L)$ and $T(R)$ when admissible path factorizes on paths for two pairs of trees of the form:

Fig.11. Factorization of paths.

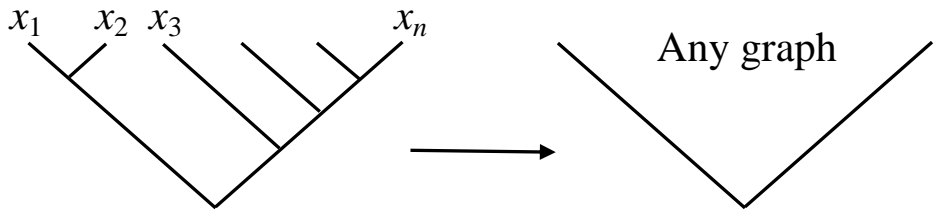

Fig.12. The second example when existence of admissible path is proved.

The corresponding path is shown in fig.12. Considering our remark we shall treat only the cases when no diagonals on polygonal representations of two trees coincide. Then there are two cases to be solved. The first one is when $(x_m\ x_{m+1})$ with $m \geq 3$ is encountered in $R$. Then the first step of transforming $T(L)$ to $T(R)$ will be organizing sub-tree for $(x_m\ x_{m+1})$ by corresponding transplantation (see fig.13a). The coloring of vortices at this step is arbitrary. Now brackets in $(x_m\ x_{m+1})$ will never be moved. Introducing $y = (x_m\ x_{m+1})$ one gets the same problem for $n - 1$ factors when one (arbitrary vortex of corresponding tree gets arbitrary sign. Note that for $n = 5$ the desired admissible path exists (base case of induction). So induction works in this case.

The second case is when $(x_2 x_3)$ is encountered in $R$ and no other $x_i$ (with $i \leq m$) are in common brackets. Then $T(R)$ is of the form shown in fig. 13b.

Now if $m < n$ use algorithm of fig.13b.

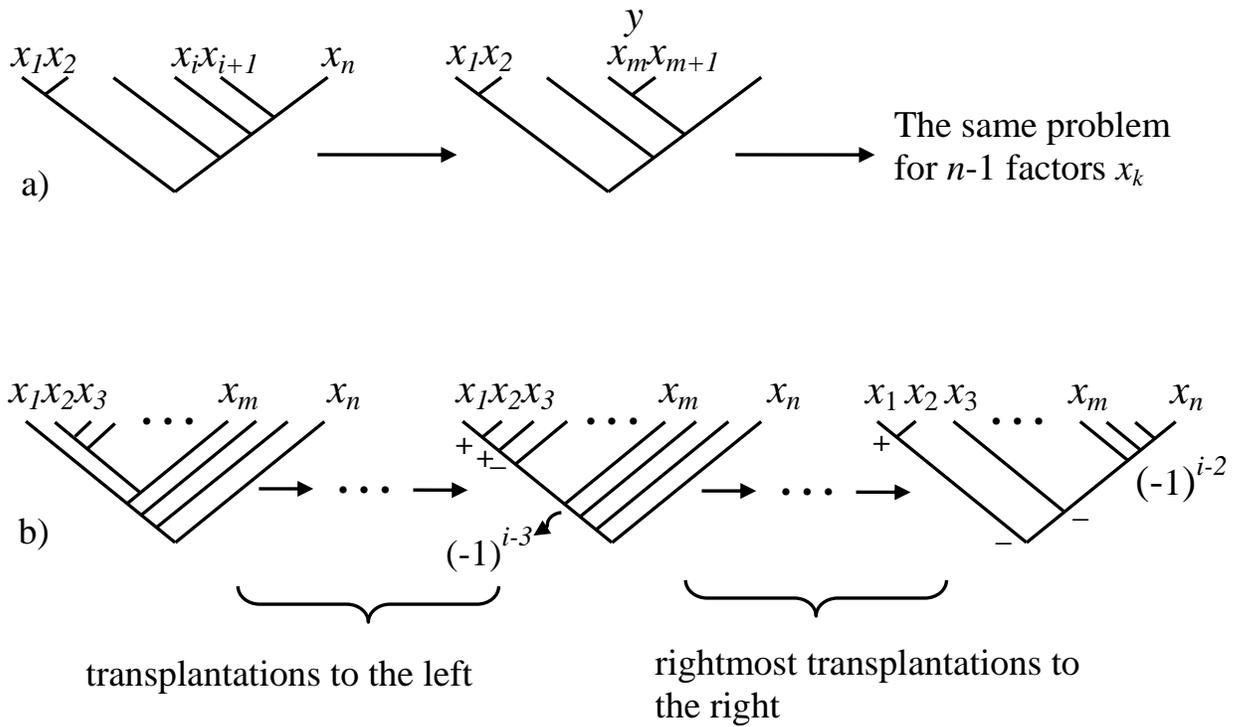

Fig. 13. Algorithm for the second example.

. If $m = n$ just do all rightmost transplantations to the right and then move the brunch with $x_2$ to the left. It is possible to do the same analysis when one of the trees in the map is like the one shown in fig.14. Corresponding polygon is also shown in fig.14. Rotating (and getting the mirror image of) this polygon one gets the trees for which the same algorithm is applicable.

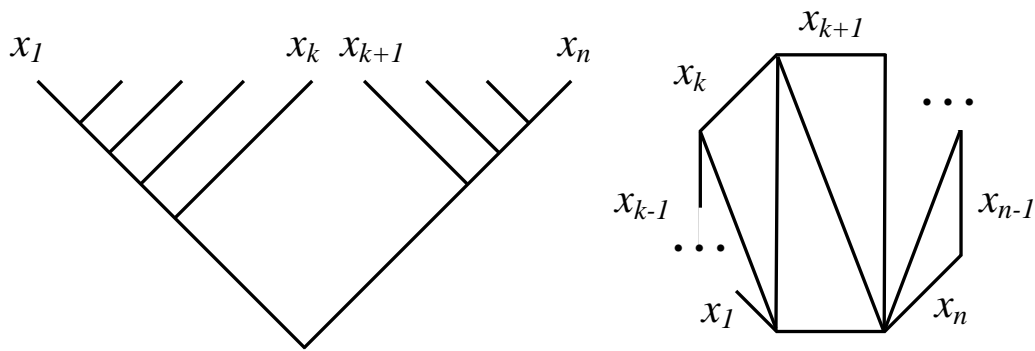

Fig.14. The class of trees in which the tree of the second example is present.

Some other examples of admissible paths give us the hope that the conjecture of this work is true.

I'd like to thank Dr. E.L. Surkov for his attention to this work.